\documentclass[12pt,reqno]{amsart}
\usepackage{amsmath, amsthm, amscd, amsfonts, amssymb, graphicx, color}
\usepackage[bookmarksnumbered, colorlinks, plainpages]{hyperref}

\theoremstyle{definition}

\theoremstyle{remark}

\numberwithin{equation}{section}
\begin{document}
\begin{center}
{\textbf{Some results on $\eta$-Yamabe Solitons in 3-dimensional trans-Sasakian manifold}}
\end{center}
\vskip 0.3cm
\begin{center}By\end{center}\vskip 0.3cm
\begin{center}
{Soumendu Roy \footnote{The first author is the corresponding author, supported by Swami Vivekananda Merit Cum Means Scholarship, Government of West Bengal, India.}, Santu Dey $^2$ and Arindam~~Bhattacharyya $^3$}
\end{center}
\vskip 0.3cm
\address[Soumendu Roy]{Department of Mathematics, Jadavpur University, Kolkata-700032, India}
\email{soumendu1103mtma@gmail.com}

\address[Santu Dey]{Department of Mathematics, Bidhan Chandra College, Asansol - 4, West Bengal-713303 , India}
\email{santu.mathju@gmail.com}

\address[Arindam Bhattacharyya]{Department of Mathematics, Jadavpur University, Kolkata-700032, India}
\email{bhattachar1968@yahoo.co.in}
\vskip 0.5cm
\begin{center}
\textbf{Abstract}\end{center}
The object of the present paper is to study some properties of 3-dimensional trans-Sasakian manifold whose metric is $\eta$-Yamabe soliton. We have studied here some certain curvature conditions of 3-dimensional trans-Sasakian manifold admitting $\eta$-Yamabe soliton. Lastly we construct a 3-dimensional trans-Sasakian manifold satisfying $\eta$-Yamabe soliton.
\vspace {0.3cm}
\section{\textbf{Introduction}}
The concept of Yamabe flow was first introduced by Hamilton \cite{hamil} to construct Yamabe metrics on compact Riemannian manifolds. On a Riemannian or pseudo-Riemannian manifold $M$, a time-dependent metric $g(\cdot, t)$ is said to evolve by the Yamabe flow if the metric $g$ satisfies the given equation,
\begin{equation}\label{1.1}
  \frac{\partial }{\partial t}g(t)=-rg(t),\hspace{0.5cm} g(0)=g_{0},
\end{equation}
where $r$ is the scalar curvature of the manifold $M$.\\\\
In 2-dimension the Yamabe flow is equivalent to the Ricci flow (defined by $\frac{\partial }{\partial t}g(t) = -2S(g(t))$, where $S$ denotes the Ricci tensor). But in dimension $> 2$ the Yamabe and Ricci flows do not agree, since the Yamabe flow preserves the conformal class of the metric but the Ricci flow does not in general.\\\\
A Yamabe soliton \cite{barbosa}, \cite{roy2}  correspond to self-similar solution of the Yamabe flow, is defined on a Riemannian or pseudo-Riemannian manifold $(M, g)$ by a vector field $\xi$ satisfying the equation,
\begin{equation}\label{1.2}
  \frac{1}{2}\pounds_V g = (r-\lambda)g,
\end{equation}
where $\pounds_\xi g$ denotes the Lie derivative of the metric $g$ along the vector field $\xi$, $r$ is the scalar curvature and $\lambda$ is a constant. Moreover a Yamabe soliton is said to be expanding if $\lambda > 0$, steady if $\lambda = 0$ and shrinking if $\lambda < 0$.\\
Yamabe solitons on a three-dimensional Sasakian manifold were studied by R. Sharma \cite{sharma}.\\
If the potential vector field $V$ is of gradient type, $V = grad(f)$, for $f$ a smooth function on $M$, then $(V, \lambda)$ is called a gradient Yamabe soliton.\\\\
\textbf{Definition 1.1:} As a generalization of Yamabe soliton, a Riemannian metric on $(M,g)$ is said to be a $\eta$-Yamabe soliton \cite{deshmukh} if
\begin{equation}\label{1.3}
  \frac{1}{2}\pounds_\xi g = (r-\lambda)g-\mu \eta \otimes \eta,
\end{equation}
where $\lambda$ and $\mu$ are contants and $\eta$ is a 1-form.\\
If $\lambda$ and $\mu$ are two smooth functions then \eqref{1.3} is said to be an almost $\eta$-Yamabe soliton or a quasi-Yamabe soliton \cite{deshmukh}.\\
Moreover if $\mu = 0$, the above equation \eqref{1.3} reduces to \eqref{1.2} and so the $\eta$-Yamabe soliton becomes Yamabe soliton. Similarly an almost $\eta$-Yamabe soliton reduces to almost Yamabe soliton if in \eqref{1.3}, $\lambda$ is a smooth function and $\mu = 0$.\\\\
Again,
\begin{equation}\label{1.4}
  R(X, Y )Z = \nabla_X\nabla_Y Z - \nabla_Y \nabla_X Z - \nabla_{[X,Y]} Z,
\end{equation}
\begin{eqnarray}\label{1.5}
 H(X, Y )Z &=& R(X, Y )Z -\frac{1}{(n-2)}[g(Y,Z)QX - g(X,Z)QY\nonumber \\
           &+& S(Y,Z)X - S(X,Z)Y],
\end{eqnarray}
 \begin{equation}\label{1.6}
       P(X, Y )Z = R(X, Y )Z -\frac{1}{(n-1)}[g(QY,Z)X - g(QX,Z)Y ],
     \end{equation}
 \begin{equation}\label{1.7}
   \tilde{C}(X,Y)Z=R(X,Y)Z-\frac{r}{n(n-1)}[g(Y,Z)X-g(X,Z)Y],
 \end{equation}
 \begin{eqnarray}\label{1.8}
   C^*(X,Y)Z &=& aR(X,Y)Z+b[S(Y,Z)X-S(X,Z)Y+g(Y,Z)QX \nonumber \\
             &-& g(X,Z)QY]-\frac{r}{n}[\frac{a}{n-1}+2b][g(Y,Z)X-g(X,Z)Y],
 \end{eqnarray}
 where $a$, $b$ are constants,
 \begin{equation}\label{1.9}
   W_2(X,Y)Z=R(X,Y)Z+\frac{1}{n-1}[g(X,Z)QY-g(Y,Z)QX],
 \end{equation}
 are the Riemannian-Christoffel curvature tensor $R$ \cite{ozgu}, the conharmonic
curvature tensor $H$ \cite{ishii}, the projective curvature tensor $P$ \cite{yano}, the concircular curvature tensor $\tilde{C}$ \cite{MISHRA}, the quasi-conformal curvature tensor $C^*$ \cite{yanosawaki} and the $W_2$-curvature tensor \cite{MISHRA} respectively in a Riemannian manifold $(M^n, g)$, where $Q$ is the Ricci operator, defined by $S(X, Y) = g(QX, Y )$, $S$ is the Ricci tensor, $r = tr(S)$ is the scalar curvature, where $tr(S)$ is the trace of $S$ and $X, Y, Z \in \chi(M)$, $\chi(M)$ being the Lie algebra of vector fields of M.\\\\
Now in \eqref{1.8} if $a = 1$ and $b = -\frac{1}{n-2}$, then we get,
\begin{eqnarray}\label{1.10}
  C^*(X,Y)Z &=& R(X,Y)Z-\frac{1}{n-2}[S(Y,Z)X-S(X,Z)Y+g(Y,Z)QX \nonumber \\
             &-& g(X,Z)QY]+\frac{r}{(n-1)(n-2)}[g(Y,Z)X-g(X,Z)Y] \nonumber \\
             &=& C(X,Y)Z,
\end{eqnarray}
where $C$ is the conformal curvature tensor \cite{Eisenhart}. Thus the conformal curvature tensor $C$ is a particular case of the tensor $C^*$.\\\\
In the present paper we study $\eta$-Yamabe soliton on 3-dimensional trans-Sasakian manifolds. The paper is organized as follows:\\
After introduction, section 2 is devoted for preliminaries on 3-dimensional trans-Sasakian manifolds. In section 3, we have studied $\eta$-Yamabe soliton on 3-dimensio-nal trans-Sasakian manifolds. Here we examine if a 3-dimensional trans-Sasakian manifold admits $\eta$-Yamabe soliton, then the scalar curvature is constant and the manifold becomes $\eta$-Einstein. We also characterized the nature of the manifold if the manifold is Ricci symmetric and the Ricci tensor is $\eta$-recurrent. Section 4 deals with the curvature properties of 3-dimensional trans-Sasakian manifold. In this section we have shown the nature of the  $\eta$-Yamabe soliton when the manifold is $\xi$- projectively flat, $\xi$-concircularly flat, $\xi$-conharmonically flat, $\xi$-quasi-conformally flat. Here we have obtained some results on $\eta$-Yamabe soliton satisfying the conditions $R(\xi,X)\cdot S = 0$ and $W_2(\xi,X)\cdot S= 0$.\\
In last section we gave an example of a 3-dimensional trans-Sasakian manifold satisfying $\eta$-Yamabe soliton.
\vspace {0.3cm}
\section{\textbf{Preliminaries}}
Let $M$ be a connected almost contact metric manifold with an almost contact metric structure $(\phi, \xi, \eta, g)$ where $\phi$ is a $(1,1)$ tensor field, $\xi$ is a vector field, $\eta$ is a 1-form  and $g$ is the compatible Riemannian metric such that
\begin{equation}\label{2.1}
\phi^2(X) = -X + \eta(X)\xi, \eta(\xi) = 1, \eta \circ \phi = 0, \phi \xi = 0,
\end{equation}
\begin{equation}\label{2.2}
g(\phi X,\phi Y) = g(X,Y) - \eta(X)\eta(Y),
\end{equation}
\begin{equation}\label{2.3}
g(X,\phi Y) = -g(\phi X,Y),
\end{equation}
\begin{equation}\label{2.4}
g(X,\xi) = \eta(X),
\end{equation}
for all vector fields $X, Y \in \chi(M).$\\
An almost contact metric structure $(\phi,\xi,\eta, g)$ on $M$ is called a trans-Sasakian structure
\cite{oja}, if $(M \times R, J,G)$ belongs to the class $W_4$ \cite{agmlh}, where $J$ is the almost complex
structure on $M \times R$ defined by $J(X, f \frac{d}{dt}) = (\phi X - f\xi, \eta(X)\frac{d}{dt})$ for all vector fields $X$ on $M$ and smooth functions $f$ on $M \times R$. It can be expressed by the condition \cite{edb},
\begin{equation}\label{2.5}
(\nabla_X \phi)Y = \alpha(g(X,Y)\xi - \eta(Y)X) + \beta(g(\phi X,Y)\xi - \eta(Y)\phi X),
\end{equation}
for some smooth functions $\alpha,\beta$ on $M$ and we say that the trans-Sasakian structure is of type $(\alpha,\beta)$.
From the above expression we can write
\begin{equation}\label{2.6}
\nabla_X \xi = -\alpha\phi X + \beta(X - \eta(X)\xi),
\end{equation}
\begin{equation}\label{2.7}
(\nabla_X \eta)Y = -\alpha g(\phi X,Y) + \beta g(\phi X,\phi Y).
\end{equation}

For a 3-dimensional trans-Sasakian manifold the following relations hold \cite{tama}, \cite{soumendu}:
\begin{equation}\label{2.8}
2\alpha \beta + \xi \alpha = 0,
\end{equation}
\begin{equation}\label{2.9}
S(X,\xi) = (2(\alpha^2 - \beta^2) - \xi \beta)\eta(X) - X\beta - (\phi X)\alpha,
\end{equation}
\begin{eqnarray}\label{2.10}
S(X,Y) &=& (\frac{r}{2} + \xi \beta - (\alpha^2 - \beta^2))g(X,Y) - (\frac{r}{2} + \xi \beta - 3(\alpha^2 - \beta^2))\eta(X)\eta(Y)\nonumber\\
&-&(Y\beta + (\phi Y)\alpha)\eta(X) - (X\beta + (\phi X)\alpha)\eta(Y),
\end{eqnarray}
where $S$ denotes the Ricci tensor of type $(0,2)$, $r$ is the scalar curvature of\\ the manifold $M$ and $\alpha, \beta$ are defined as earlier.\\
For $\alpha, \beta$ = constant the following relations hold \cite{tama}, \cite{soumendu}:
\begin{equation}\label{2.11}
S(X,Y) = (\frac{r}{2} - (\alpha^2 - \beta^2))g(X,Y) - (\frac{r}{2} - 3(\alpha^2 - \beta^2))\eta(X)\eta(Y),
\end{equation}
\begin{equation}\label{2.12}
S(X,\xi) = 2(\alpha^2 - \beta^2)\eta(X),
\end{equation}
\begin{equation}\label{2.13}
R(X,Y)\xi = (\alpha^2 - \beta^2)[\eta(Y)X - \eta(X)Y],
\end{equation}
\begin{equation}\label{2.14}
R(\xi,X)Y=  (\alpha^2 - \beta^2) [g(X,Y)\xi- \eta(Y)X],
\end{equation}
\begin{equation}\label{2.15}
R(\xi,X)\xi= (\alpha^2 - \beta^2) [\eta(X)\xi- X],
\end{equation}
\begin{equation}\label{2.16}
\eta(R(X,Y)Z)= (\alpha^2 - \beta^2) [g(Y,Z)\eta(X)- g(X,Z)\eta(Y)],
\end{equation}
where $R$ is the Riemannian curvature tensor.
\begin{equation}\label{2.17}
QX = (\frac{r}{2} - (\alpha^2 - \beta^2))X - (\frac{r}{2} - 3(\alpha^2 - \beta^2))\eta(X)\xi,
\end{equation}
where $Q$ is the Ricci operator defined earlier.\\
Again,
\begin{eqnarray}
(\pounds_\xi g)(X,Y) &= &(\nabla_\xi g)( X,Y) -\alpha g(\phi X,Y) + 2\beta g(X,Y) - 2\beta \eta(X)\eta(Y) \nonumber\\
&-& \alpha g(X,\phi Y). \nonumber
\end{eqnarray}
Then using \eqref{2.3}, the above equation becomes,
\begin{equation}\label{2.18}
  (\pounds_\xi g)(X,Y)=2\beta g(X,Y) - 2\beta \eta(X)\eta(Y).
\end{equation}
where $\nabla$ is the Levi-Civita  connection associated with $g$ and $\pounds_\xi$ denotes the Lie derivative along the vector field $\xi$.
\vspace {0.3cm}
\section{\textbf{$\eta$-Yamabe soliton on 3-dimensional trans-Sasakian manifold}}
Let $M$ be a 3-dimensional trans-Sasakian manifold. Consider the $\eta$-Yamabe soliton on $M$ as:
\begin{equation}\label{3.1}
\frac{1}{2}(\pounds_\xi g)(X,Y) = (r-\lambda)g(X,Y)-\mu \eta(X) \eta(Y),
\end{equation}
for all vector fields $X, Y$ on $M$.\\
Then from \eqref{2.18} and \eqref{3.1}, we get,
\begin{equation}\label{3.2}
 (r-\lambda-\beta)g(X,Y)=(\mu-\beta)\eta(X)\eta(Y).
\end{equation}
Taking $Y=\xi$ in the above equation and using \eqref{2.1}, we have,
\begin{equation}\label{3.3}
  (r-\lambda-\mu)\eta(X)=0.
\end{equation}
Since $\eta(X) \neq 0$, so we get,
\begin{equation}\label{3.4}
  r=\lambda+\mu.
\end{equation}
Now as both $\lambda$ and $\mu$ are constants, $r$ is also constant.\\\\
So we can state the following theorem:\\
\textbf{Theorem 3.1.} {\em If a 3-dimensional trans-Sasakian manifold $M$ admits an $\eta$-Yamabe soliton $(g, \xi)$, $\xi$ being the Reeb vector field of $M$, then the scalar curvature is constant.}\\\\
In \eqref{3.4} if $\mu = 0$, we get $r = \lambda$ and so \eqref{3.1} becomes,
\begin{equation}
  \pounds_\xi g=0. \nonumber
\end{equation}
Thus $\xi$ is a Killing vector field and consequently $M$ is a 3-dimensional $K$-trans-Sasakian manifold.\\\\
Then we have,\\
\textbf{Corollary 3.2.} {\em If a 3-dimensional trans-Sasakian manifold $M$ admits a Yamabe soliton $(g, \xi)$, $\xi$ being the Reeb vector field of $M$, then the manifold is a 3-dimensional $K$-trans-Sasakian manifold.}\\\\
Now from \eqref{2.11} and \eqref{3.4}, we have,
\begin{equation}\label{3.5}
  S(X,Y) = (\frac{\lambda+\mu}{2} - (\alpha^2 - \beta^2))g(X,Y) - (\frac{\lambda+\mu}{2} - 3(\alpha^2 - \beta^2))\eta(X)\eta(Y),
\end{equation}
for all vector fields $X, Y$ on $M$.\\\\
This leads to the following:\\
\textbf{Corollary 3.3.} {\em If a 3-dimensional trans-Sasakian manifold $M$ admits an $\eta$-Yamabe soliton $(g, \xi)$, $\xi$ being the Reeb vector field of $M$, then the manifold becomes $\eta$-Einstein manifold.}\\\\
We know,
\begin{equation}\label{3.6}
  (\nabla_X S)(Y,Z)= XS(Y,Z)-S(\nabla_X Y,Z)-S(Y,\nabla_X Z),
\end{equation}
for all vector fields $X, Y, Z$ on $M$ and $\nabla$ is the Levi-Civita  connection associated with $g$.\\
Now replacing the expression of S from \eqref{3.5}, we obtain,
\begin{equation}\label{3.7}
  (\nabla_X S)(Y,Z) = -[\frac{\lambda+\mu}{2} - 3(\alpha^2 - \beta^2)][\eta(Z)(\nabla_X \eta)Y+\eta(Y)(\nabla_X \eta)Z].
\end{equation}
for all vector fields $X, Y, Z$ on $M$.\\\\
Now let the manifold be Ricci symmetric i.e $ \nabla S = 0.$\\
Then from \eqref{3.7}, we have,
\begin{equation}\label{3.8}
  [\frac{\lambda+\mu}{2} - 3(\alpha^2 - \beta^2)][\eta(Z)(\nabla_X \eta)Y+\eta(Y)(\nabla_X \eta)Z]=0,
\end{equation}
for all vector fields $X, Y, Z$ on $M$.\\
Taking $Z = \xi$ in the above equation and using \eqref{2.7}, \eqref{2.1}, we get,
\begin{equation}\label{3.9}
  [\frac{\lambda+\mu}{2} - 3(\alpha^2 - \beta^2)][\beta g(\phi X,\phi Y)-\alpha g(\phi X, Y)]=0
\end{equation}
for all vector fields $X, Y$ on $M$.\\
Hence we get,
\begin{equation}
  \lambda+\mu=6(\alpha^2 - \beta^2). \nonumber
\end{equation}
This leads to the following:\\
\textbf{Proposition 3.4.} {\em Let a 3-dimensional trans-Sasakian manifold $M$ admit an $\eta$-Yamabe soliton $(g, \xi)$, $\xi$ being the Reeb vector field of $M$. If the manifold is Ricci symmetric then $\lambda+\mu=6(\alpha^2 - \beta^2)$, where $\lambda, \mu, \alpha, \beta$ are constants.}\\\\
Now if the Ricci tensor $S$ is $\eta$-recurrent, then we have,
\begin{equation}\label{3.10}
  \nabla S= \eta \otimes S,
\end{equation}
which implies that,
\begin{equation}\label{3.11}
  (\nabla_X S)(Y,Z)=\eta(X)S(Y,Z),
\end{equation}
for all vector fields $X, Y, Z$ on $M$.\\
Then using \eqref{3.7}, we get,
\begin{equation}\label{3.12}
  -[\frac{\lambda+\mu}{2} - 3(\alpha^2 - \beta^2)][\eta(Z)(\nabla_X \eta)Y+\eta(Y)(\nabla_X \eta)Z]=\eta(X)S(Y,Z),
\end{equation}
for all vector fields $X, Y, Z$ on $M$.\\
Using \eqref{2.7}, the above equation becomes,
\begin{eqnarray}\label{3.13}
  &&-[\frac{\lambda+\mu}{2} - 3(\alpha^2 - \beta^2)][\eta(Z)(-\alpha g(\phi X,Y) + \beta g(\phi X,\phi Y)) \nonumber\\
  &&+\eta(Y)(-\alpha g(\phi X,Z)+\beta g(\phi X,\phi Z))]=\eta(X)S(Y,Z).
\end{eqnarray}
Now taking $Y = \xi, Z = \xi$ and using \eqref{2.1}, \eqref{3.5}, the above equation becomes,
\begin{equation}
  2(\alpha^2-\beta^2)\eta(X)=0.\nonumber
\end{equation}
Since $\eta(X)  \neq 0$, for all $X$ on $M$, we have,
\begin{equation}\label{3.14}
  \alpha=\pm \beta.
\end{equation}
This leads to the following:\\
\textbf{Proposition 3.5.} {\em  Let a 3-dimensional trans-Sasakian manifold $M$ admit an $\eta$-Yamabe soliton $(g, \xi)$, $\xi$ being the Reeb vector field of $M$. If the Ricci tensor $S$ is $\eta$-recurrent then $\alpha=\pm \beta$.}\\\\
Now if the manifold is Ricci symmetric and the Ricci tensor $S$ is $\eta$-recurrent, then using \eqref{3.14} in $\lambda+\mu=6(\alpha^2 - \beta^2)$ and from \eqref{3.4}, we have the following:\\\\
\textbf{Proposition 3.6.} {Let a 3-dimensional trans-Sasakian manifold $M$ admit an $\eta$-Yamabe soliton $(g, \xi)$, $\xi$ being the Reeb vector field of $M$. If the manifold is Ricci symmetric and the Ricci tensor $S$ is $\eta$-recurrent then the manifold becomes flat.}\\\\
Let an $\eta$-Yamabe soliton be defined on a 3-dimensional trans-Sasakian manifold $M$ as,
\begin{equation}\label{3.15}
  \frac{1}{2}\pounds_V g = (r-\lambda)g-\mu \eta \otimes \eta,
\end{equation}
where $\pounds_V g$ denotes the Lie derivative of the metric $g$ along a vector field $V$, $r$ is defined as \eqref{1.2} and $\lambda, \mu$ are defined as \eqref{1.3}.\\
Let $V$ be pointwise co-linear with $\xi$ i.e. $V = b\xi$ where $b$ is a function on $M$.\\
Then the equation \eqref{3.15} becomes,
\begin{equation}\label{3.16}
  (\pounds_{b\xi} g)(X,Y) = 2(r-\lambda)g(X,Y)-2\mu \eta(X) \eta(Y),,
\end{equation}
for any vector fields $X, Y$ on $M$.\\
Applying the property of Lie derivative and Levi-Civita connection we have,
\begin{equation}\label{3.17}
  bg(\nabla_X \xi,Y)+(Xb)\eta(Y)+bg(\nabla_Y \xi,X)+(Yb)\eta(X)= 2(r-\lambda)g(X,Y)-2\mu \eta(X) \eta(Y).
\end{equation}
Using \eqref{2.6} and \eqref{2.3}, the above equation reduces to,
\begin{equation}\label{3.18}
  2b\beta [g(X,Y)- \eta(X) \eta(Y)]+(Xb)\eta(Y)+(Yb)\eta(X)= 2(r-\lambda)g(X,Y)-2\mu \eta(X) \eta(Y).
\end{equation}
Now taking $Y = \xi$ in the above equation and using \eqref{2.1}, \eqref{2.4}, we obtain,
\begin{equation}\label{3.19}
  Xb+(\xi b)\eta(X)=2(r-\lambda)\eta(X)-2\mu \eta(X).
\end{equation}
Again taking $X = \xi$, we get,
\begin{equation}\label{3.20}
  \xi b=r-\lambda-\mu.
\end{equation}
Then using \eqref{3.20}, the equation \eqref{3.19} becomes,
\begin{equation}\label{3.21}
  Xb=(r-\lambda-\mu)\eta(X).
\end{equation}
Applying exterior differentiation in \eqref{3.21}, we have,
\begin{equation}\label{3.22}
  (r-\lambda-\mu)d \eta=0.
\end{equation}
Since $d\eta \neq 0$ \cite{tama}, the above equation gives,
\begin{equation}\label{3.23}
  r=\lambda+\mu.
\end{equation}
Using \eqref{3.23}, the equation \eqref{3.21} becomes,
\begin{equation}\label{3.24}
  Xb=0,
\end{equation}
which implies that $b$ is constant.\\\\
So we can state the following theorem:\\
\textbf{Theorem 3.7.} {\em Let $M$ be a 3-dimensional trans-Sasakian manifold admitting an $\eta$-Yamabe soliton $(g, V)$, $V$ being a vector field on $M$. If $V$ is pointwise co-linear with $\xi$, then $V$ is a constant multiple of $\xi$, where $\xi$ being the Reeb vector field of $M$.}\\\\
Using \eqref{3.23}, the equation \eqref{3.15} becomes,
\begin{equation}\label{3.25}
  (\pounds_V g)(X,Y)=2\mu[g(X,Y)-\eta(X)\eta(Y)],
\end{equation}
for all vector fields $X, Y, Z$ on $M$.\\\\
Then we have,\\
\textbf{Corollary 3.8.} {\em Let $M$ be a 3-dimensional trans-Sasakian manifold admitting an $\eta$-Yamabe soliton $(g, V)$, $V$ being a vector field on $M$ which is pointwise co-linear with $\xi$, where $\xi$ being the Reeb vector field of $M$. V is a Killing vector field iff the soliton reduces to a Yamabe soliton.}\\\\
From the equation \eqref{3.5}, we get,
\begin{equation}\label{3.26}
  QX=(\frac{\lambda+\mu}{2} - (\alpha^2 - \beta^2))X - (\frac{\lambda+\mu}{2} - 3(\alpha^2 - \beta^2))\eta(X)\xi,
\end{equation}
for any vector field $X$ on $M$ and $Q$ is defined as earlier.\\
We know,
 \begin{equation}\label{3.27}
    (\nabla_\xi Q)X = \nabla_\xi QX - Q(\nabla_\xi X),
 \end{equation}
for any vector field $X$ on $M$.\\
Then using \eqref{3.26}, the equation \eqref{3.27} becomes,
\begin{equation}\label{3.28}
  (\nabla_\xi Q)X=-[\frac{\lambda+\mu}{2}-3(\alpha^2-\beta^2)]((\nabla_\xi \eta)X)\xi.
\end{equation}
Using \eqref{2.7} in the above equation, we get,
\begin{equation}\label{3.29}
  (\nabla_\xi Q)X=0,
\end{equation}
for any vector field $X$ on $M$.\\
Hence $Q$ is parallel along $\xi$.\\
Again from \eqref{3.7}, we obtain,
\begin{equation}\label{3.30}
  (\nabla_\xi S)(X,Y) = -[\frac{\lambda+\mu}{2} - 3(\alpha^2 - \beta^2)][\eta(Y)(\nabla_\xi \eta)X+\eta(X)(\nabla_\xi \eta)Y],
\end{equation}
for any vector fields $X, Y$ on $M$.\\
Using \eqref{2.7} in the above equation, we get,
\begin{equation}\label{3.31}
  (\nabla_\xi S)(X,Y)=0,
\end{equation}
for any vector fields $X, Y$ on $M$.\\
Hence $S$ is parallel along $\xi$.\\\\
So we can state the following theorem:\\
\textbf{Theorem 3.9.} {\em Let $M$ be a 3-dimensional trans-Sasakian manifold admitting an $\eta$-Yamabe soliton $(g, \xi)$, $\xi$ being the Reeb vector field on $M$. Then $Q$ and $S$ are parallel along $\xi$, where $Q$ is the Ricci operator, defined by $S(X, Y) = g(QX, Y )$ and $S$ is the Ricci tensor of $M$.}
\vspace {0.3cm}
\section{\textbf{Curvature properties on 3-dimensional trans-Sasakian manifold  admitting $\eta$-Yamabe soliton}}
From the definition of projective curvature tensor \eqref{1.6}, defined on a 3-dimensio-nal trans-Sasakian manifold and using the property $g(QX, Y) = S(X, Y)$, we have,
\begin{equation}\label{4.1}
  P(X, Y )Z = R(X, Y )Z -\frac{1}{2}[S(Y,Z)X - S(X,Z)Y ],
\end{equation}
for any vector fields $X, Y, Z$ on $M$.\\
Putting $Z = \xi$ in the above equation and using \eqref{2.13} and \eqref{3.5}, we obtain,
\begin{equation}\label{4.2}
  P(X, Y)\xi=(\alpha^2 - \beta^2)[\eta(Y)X - \eta(X)Y]-\frac{1}{2}[2(\alpha^2-\beta^2)\eta(Y)X-2(\alpha^2-\beta^2)\eta(X)Y],
\end{equation}
which implies that,
\begin{equation}\label{4.3}
  P(X, Y)\xi=0.
\end{equation}
So we can state the following theorem:\\
\textbf{Theorem 4.1.} {\em A 3-dimensional trans-Sasakian manifold $M$  admitting $\eta$-Yamabe soliton $(g, \xi)$, $\xi$ being the Reeb vector field on $M$ is $\xi$- projectively flat.}\\\\
From the definition of concircular curvature tensor \eqref{1.7}, defined on a 3-dimensio-nal trans-Sasakian manifold, we have,
\begin{equation}\label{4.4}
  \tilde{C}(X,Y)Z=R(X,Y)Z-\frac{r}{6}[g(Y,Z)X-g(X,Z)Y],
\end{equation}
for any vector fields $X, Y, Z$ on $M$.\\
Putting $Z = \xi$ in the above equation and using \eqref{2.4} and \eqref{2.13}, we obtain,
\begin{equation}\label{4.5}
  \tilde{C}(X,Y)\xi=(\alpha^2 - \beta^2)[\eta(Y)X - \eta(X)Y]-\frac{r}{6}[\eta(Y)X - \eta(X)Y].
\end{equation}
Now using \eqref{3.4}, we get,
\begin{equation}\label{4.6}
  \tilde{C}(X,Y)\xi=[(\alpha^2 - \beta^2)-\frac{\lambda+\mu}{6}][\eta(Y)X - \eta(X)Y].
\end{equation}
This implies that $\tilde{C}(X, Y )\xi = 0$ iff $\lambda+\mu = 6(\alpha^2 - \beta^2).$\\\\
So we can state the following theorem:\\
\textbf{Theorem 4.2.} {\em A 3-dimensional trans-Sasakian manifold $M$  admitting $\eta$-Yamabe soliton $(g, \xi)$, $\xi$ being the Reeb vector field on $M$ is $\xi$-concircularly flat iff $\lambda+\mu = 6(\alpha^2 - \beta^2).$}\\\\
Now if the Ricci tensor $S$ is $\eta$- recurrent then using \eqref{3.14} in \eqref{4.5}, we have,\\\\
\textbf{Corollary 4.3.} {\em Let $M$ be a 3-dimensional trans-Sasakian manifold admitting an $\eta$-Yamabe soliton $(g, \xi)$, $\xi$ being the Reeb vector field on $M$. If the manifold is $\xi$-concircularly flat and the Ricci tensor is $\eta$- recurrent then the manifold $M$ becomes flat.}\\\\
From the definition of conharmonic curvature tensor \eqref{1.5}, defined on a 3-dimensional trans-Sasakian manifold, we have,
\begin{equation}\label{4.7}
  H(X, Y )Z = R(X, Y )Z -[g(Y,Z)QX - g(X,Z)QY + S(Y,Z)X - S(X,Z)Y],
\end{equation}
for any vector fields $X, Y, Z$ on $M$.\\
Putting $Z = \xi$ in the above equation and using \eqref{2.4}, \eqref{2.13}, \eqref{3.5} and \eqref{3.26}, the above equation becomes,
\begin{equation}\label{4.8}
  H(X, Y )\xi =(\alpha^2 - \beta^2)[\eta(Y)X - \eta(X)Y]-[\frac{\lambda+\mu}{2}+(\alpha^2-\beta^2)][\eta(Y)X - \eta(X)Y].
\end{equation}
Hence we get,
\begin{equation}\label{4.9}
  H(X, Y )\xi =-\frac{\lambda+\mu}{2}[\eta(Y)X - \eta(X)Y].
\end{equation}
This implies that $H(X, Y )\xi = 0$ iff $\lambda+\mu = 0.$\\\\
So we can state the following theorem:\\
\textbf{Theorem 4.4.} {\em A 3-dimensional trans-Sasakian manifold $M$  admitting $\eta$-Yamabe soliton $(g, \xi)$, $\xi$ being the Reeb vector field on $M$ is $\xi$-conharmonically flat iff $\lambda+\mu = 0.$}\\\\
From the definition of quasi-conformal curvature tensor \eqref{1.8}, defined on a 3-dimensional trans-Sasakian manifold, we have,
\begin{eqnarray}\label{4.10}
   C^*(X,Y)Z &=& aR(X,Y)Z+b[S(Y,Z)X-S(X,Z)Y+g(Y,Z)QX \nonumber \\
             &-& g(X,Z)QY]-\frac{r}{3}[\frac{a}{2}+2b][g(Y,Z)X-g(X,Z)Y],
 \end{eqnarray}
for any vector fields $X, Y, Z$ on $M$ and $a$, $b$ are constants.\\
Putting $Z = \xi$ in the above equation and using \eqref{2.4}, \eqref{2.13}, \eqref{3.4}, \eqref{3.5} and \eqref{3.26}, the above equation becomes,
\begin{eqnarray}\label{4.11}
   C^*(X,Y)\xi &=& a(\alpha^2 - \beta^2)[\eta(Y)X - \eta(X)Y] \nonumber \\
               &+& b[\frac{\lambda+\mu}{2}+(\alpha^2-\beta^2)][\eta(Y)X - \eta(X)Y] \nonumber\\
               &-& \frac{\lambda+\mu}{3}[\frac{a}{2}+2b][\eta(Y)X - \eta(X)Y].
 \end{eqnarray}
Hence we have,
\begin{equation}\label{4.12}
  C^*(X,Y)\xi = [a(\alpha^2 - \beta^2)+b[\frac{\lambda+\mu}{2}+(\alpha^2-\beta^2)]-\frac{\lambda+\mu}{3}[\frac{a}{2}+2b]][\eta(Y)X - \eta(X)Y].
\end{equation}\\
 This implies that $C^*(X,Y)\xi = 0$ iff $ a(\alpha^2 - \beta^2) + b[\frac{\lambda + \mu}{2} + (\alpha^2 - \beta^2)] - \frac{\lambda + \mu}{3}[\frac{a}{2} + 2b] = 0.$\\\\
 Then by simplifying, we obtain, $C^*(X,Y)\xi = 0$ iff $(a + b)[(\alpha^2 - \beta^2) - \frac{\lambda+\mu}{6}] = 0,$ i.e either $a + b = 0$ or $\lambda + \mu = 6(\alpha^2 - \beta^2).$\\\\
So we can state the following theorem:\\
\textbf{Theorem 4.5.} {\em A 3-dimensional trans-Sasakian manifold $M$  admitting $\eta$-Yamabe soliton $(g, \xi)$, $\xi$ being the Reeb vector field on $M$ is $\xi$-quasi-conformally flat iff either $a + b = 0$ or $\lambda + \mu = 6(\alpha^2 - \beta^2).$}\\\\
Now if the Ricci tensor $S$ is $\eta$-recurrent then using \eqref{3.14} in \eqref{4.12}, we get,
\begin{equation}\label{4.13}
  C^*(X,Y)\xi =-\frac{a+b}{6}(\lambda+\mu)[\eta(Y)X - \eta(X)Y].
\end{equation}
Hence using \eqref{3.4} in \eqref{4.13}, we have,\\\\
\textbf{Corollary 4.6.} {\em Let a 3-dimensional trans-Sasakian manifold $M$ admit an $\eta$-Yamabe soliton $(g, \xi)$, $\xi$ being the Reeb vector field on $M$. If the manifold is $\xi$-quasi-conformally flat and the Ricci tensor is $\eta$-recurrent then the manifold $M$ becomes flat, provided $a + b \neq 0.$}\\\\
We know,
\begin{equation}\label{4.14}
  R(\xi,X)\cdot S=S(R(\xi,X)Y,Z)+S(Y,R(\xi,X)Z),
\end{equation}
for any vector fields $X, Y, Z$ on $M$.\\
Now let the manifold be $\xi$-semi symmetric, i.e $R(\xi,X)\cdot S = 0.$\\
Then from \eqref{4.14}, we have,
\begin{equation}\label{4.15}
  S(R(\xi,X)Y,Z)+S(Y,R(\xi,X)Z)=0,
\end{equation}
for any vector fields $X, Y, Z$ on $M$.\\
Using \eqref{2.14}, the above equation becomes,
\begin{equation}\label{4.16}
  S((\alpha^2 - \beta^2) (g(X,Y)\xi- \eta(Y)X),Z)+S(Y,(\alpha^2 - \beta^2) (g(X,Z)\xi- \eta(Z)X))=0.
\end{equation}
Replacing the expression of S from \eqref{3.5} and simplifying we get,
\begin{equation}\label{4.17}
  (\alpha^2-\beta^2)[\frac{\lambda+\mu}{2}-3(\alpha^2-\beta^2)][g(X,Y)\eta(Z)+g(X,Z)\eta(Y)-2\eta(X)\eta(Y)\eta(Z)]=0.
\end{equation}
Taking $Z = \xi$ in the above equation and using \eqref{2.1}, \eqref{2.4}, we obtain,
\begin{equation}\label{4.18}
  (\alpha^2-\beta^2)[\frac{\lambda+\mu}{2}-3(\alpha^2-\beta^2)][g(X,Y)-\eta(x)\eta(Y)]=0,
\end{equation}
for any vector fields $X, Y$ on $M$.\\
Using \eqref{2.2}, the above equation becomes,
\begin{equation}\label{4.19}
  (\alpha^2-\beta^2)[\frac{\lambda+\mu}{2}-3(\alpha^2-\beta^2)]g(\phi X, \phi Y)=0,
\end{equation}
for any vector fields $X, Y$ on $M$.\\
Hence we get,
\begin{equation}\label{4.20}
  (\alpha^2-\beta^2)[\frac{\lambda+\mu}{2}-3(\alpha^2-\beta^2)]=0.
\end{equation}
Then either $(\alpha^2 - \beta^2) = 0$, or $\lambda + \mu = 6(\alpha^2 - \beta^2).$\\\\
So we can state the following theorem:\\
\textbf{Theorem 4.7.} {\em If a 3-dimensional trans-Sasakian manifold $M$  admitting $\eta$-Yamabe soliton $(g, \xi)$, $\xi$ being the Reeb vector field on $M$ is $\xi$-semi symmetric then either $(\alpha^2 - \beta^2) = 0$, or $\lambda + \mu = 6(\alpha^2 - \beta^2).$}\\\\
From the definition of $W_2$-curvature tensor \eqref{1.9}, defined on a 3-dimensional trans-Sasakian manifold, we have,
\begin{equation}\label{4.21}
   W_2(X,Y)Z=R(X,Y)Z+\frac{1}{2}[g(X,Z)QY-g(Y,Z)QX],
 \end{equation}
 for any vector fields $X, Y, Z$ on $M$.\\
Again we know,
\begin{equation}\label{4.22}
  W_2(\xi,X)\cdot S=S(W_2(\xi,X)Y,Z)+S(Y,W_2(\xi,X)Z),
\end{equation}
for any vector fields $X, Y, Z$ on $M$.\\
Replacing the expression of S from \eqref{3.5}, we get on simplifying,
\begin{eqnarray}\label{4.23}
   W_2(\xi,X)\cdot S &=&  [\frac{\lambda+\mu}{2}-(\alpha^2-\beta^2)][g(W_2(\xi,X)Y,Z)+g(Y,W_2(\xi,X)Z)]\nonumber\\
                     &-&  [\frac{\lambda+\mu}{2}-3(\alpha^2-\beta^2)][\eta(W_2(\xi,X)Y)\eta(Z) \nonumber \\
                     &+&   \eta(Y)\eta(W_2(\xi,X)Z)].
\end{eqnarray}
Now from the definition of $W_2$-curvature tensor \eqref{4.21} and then by using \eqref{2.14}, the property $g(QX, Y) = S(X, Y)$, \eqref{3.5}, the above equation becomes,
\begin{eqnarray}\label{4.24}
  W_2(\xi,X)\cdot S &=&\frac{1}{2}[\frac{\lambda+\mu}{2}-(\alpha^2-\beta^2)][\frac{\lambda+\mu}{2}-3(\alpha^2-\beta^2)][g(X,Y)\eta(Z) \nonumber\\
                    &+& g(X,Z)\eta(Y)-2\eta(X)\eta(Y)\eta(Z)],
\end{eqnarray}
for any vector fields $X, Y, Z$ on $M$.\\\\
Let in this manifold $M$, $W_2(\xi,X)\cdot S= 0.$\\
Then from \eqref{4.24}, we get,
\begin{eqnarray}\label{4.25}
 &&[\frac{\lambda+\mu}{2}-(\alpha^2-\beta^2)][\frac{\lambda+\mu}{2}-3(\alpha^2-\beta^2)][g(X,Y)\eta(Z) \nonumber\\
 &&+g(X,Z)\eta(Y)-2\eta(X)\eta(Y)\eta(Z)]=0,
\end{eqnarray}
for any vector fields $X, Y, Z$ on $M$.\\
Taking $Z = \xi$ in the above equation and using \eqref{2.1}, \eqref{2.4}, we obtain,
\begin{equation}\label{4.26}
   [\frac{\lambda+\mu}{2}-(\alpha^2-\beta^2)][\frac{\lambda+\mu}{2}-3(\alpha^2-\beta^2)][g(X,Y)-\eta(X)\eta(Y)]=0,
\end{equation}
for any vector fields $X, Y$ on $M$.\\
Using \eqref{2.2}, the above equation becomes,
\begin{equation}\label{4.27}
   [\frac{\lambda+\mu}{2}-(\alpha^2-\beta^2)][\frac{\lambda+\mu}{2}-3(\alpha^2-\beta^2)]g(\phi X, \phi Y)=0,
\end{equation}
for any vector fields $X, Y$ on $M$.\\
Hence we get,
\begin{equation}\label{4.28}
  [\frac{\lambda+\mu}{2}-(\alpha^2-\beta^2)][\frac{\lambda+\mu}{2}-3(\alpha^2-\beta^2)]=0.
\end{equation}
Then either $\lambda + \mu = 2(\alpha^2 - \beta^2),$ or $\lambda + \mu = 6(\alpha^2 - \beta^2).$\\\\
So we can state the following theorem:\\
\textbf{Theorem 4.8.} {\em If a 3-dimensional trans-Sasakian manifold $M$  admits an $\eta$-Yamabe soliton $(g, \xi)$, $\xi$ being the Reeb vector field on $M$ and satisfies $W_2(\xi,X)\cdot S= 0,$ where $W_2$ is the $W_2$-curvature tensor and $S$ is the Ricci tensor then either $\lambda + \mu = 2(\alpha^2 - \beta^2),$ or $\lambda + \mu = 6(\alpha^2 - \beta^2).$ }\\\\
Now if the Ricci tensor $S$ is $\eta$-recurrent then using \eqref{3.14} in \eqref{4.28} and from \eqref{3.4}, we have,\\\\
\textbf{Corollary 4.9.} {\em If a 3-dimensional trans-Sasakian manifold $M$  admits an $\eta$-Yamabe soliton $(g, \xi)$, $\xi$ being the Reeb vector field on $M$ and satisfies $W_2(\xi,X)\cdot S= 0,$ where $W_2$ is the $W_2$-curvature tensor and $S$ is the Ricci tensor which is $\eta$- recurrent, then the manifold becomes flat.}
\vspace {0.3cm}
\section{\textbf{Example of a 3-dimensional trans-Sasakian manifold admitting $\eta$-Yamabe soliton:}}
In this section we give an example of a 3-dimensional trans-Sasakian manifold with $\alpha, \beta = $constant.\\
We consider the 3-dimensional manifold $M = \{(x, y, z) \in \mathbb{R}^3, z \neq 0 \}$, where $(x, y, z)$ are standard coordinates in $\mathbb{R}^3$. Let ${e_1, e_2, e_3}$ be a linearly independent system of vector fields on $M$ given by,
\begin{equation}
   e_1=z\frac{\partial}{\partial x},\quad e_2 =z \frac{\partial}{\partial y}, \quad e_3 =z\frac{\partial}{\partial z}. \nonumber
\end{equation}
Let $g$ be the Riemannian metric defined by,
\begin{equation}\label{5.1}
  g(e_1,e_1)= g(e_2,e_2) = g(e_3,e_3)=1,\nonumber
\end{equation}
\begin{equation}\label{5.2}
   g(e_1,e_2) = g(e_2,e_3)= g(e_3,e_1) =0.\nonumber
\end{equation}
Let $\eta$ be the 1-form defined by $\eta(Z) = g(Z,e_3)$, for any $Z \in \chi(M)$,where $\chi(M)$ is the set of all differentiable vector fields on $M$ and $\phi$ be the (1, 1)-tensor field defined by,
\begin{equation}
  \phi e_1=-e_2, \quad \phi e_2=e_1,\quad  \phi e_3=0.\nonumber
\end{equation}
Then, using the linearity of $\phi$ and $g$, we have $$\eta(e_3) = 1, \phi ^2 (Z) = -Z + \eta(Z)e_3$$ and $$g(\phi Z,\phi W) = g(Z,W) - \eta(Z)\eta(W),$$ for any $Z,W \in \chi(M)$.\\
Let $\nabla$ be the Levi-Civita connection with respect to the Riemannian metric $g$. Then we have,
  $$ [e_1,e_2] = 0, [e_2,e_3] = -e_2, [e_1,e_3] = -e_1.$$
The connection $\nabla$ of the metric $g$ is given by,
\begin{eqnarray}
  2g(\nabla_X Y,Z) &=& Xg(Y,Z)+Yg(Z,X)-Zg(X,Y)\nonumber \\
                   &-& g(X, [Y,Z])-g(Y, [X, Z]) + g(Z, [X, Y]),\nonumber
\end{eqnarray}
which is known as Koszul’s formula.\\
Using Koszul’s formula, we can easily calculate,
$$\nabla_{e_1} e_3 = -e_1,  \nabla_{e_2} e_3 = -e_2 ,  \nabla_{e_3} e_3 = 0 ,$$
$$\nabla_{e_1} e_1 = e_3,  \nabla_{e_2} e_1 = 0,  \nabla_{e_3} e_1 = 0,$$
$$\nabla_{e_1} e_2 = 0,  \nabla_{e_2} e_2 = e_3,  \nabla_{e_3} e_2 = 0.$$
We see that,
\begin{multline}\label{5.1}
 (\nabla_{e_1} \phi)e_1 = \nabla_{e_1} \phi e_1-\phi \nabla_{e_1} e_1 = -\nabla_{e_1}e_2-\phi e_3= 0 \\
  =0(g(e_1,e_1)e_3 - \eta(e_1)e_1) - 1(g(\phi e_1,e_1)e_3 - \eta(e_1)\phi e_1).\quad \quad \quad \quad
\end{multline}
\begin{multline}\label{5.2}
  (\nabla_{e_1} \phi)e_2 = \nabla_{e_1} \phi e_2-\phi \nabla_{e_1} e_2 = \nabla_{e_1} e_1-0 = e_3 \\
  =0(g(e_1,e_2)e_3 - \eta(e_2)e_1) - 1(g(\phi e_1,e_2)e_3 - \eta(e_2)\phi e_1).\quad \quad \quad \quad
\end{multline}
\begin{multline}\label{5.3}
 (\nabla_{e_1} \phi)e_3 = \nabla_{e_1} \phi e_3-\phi \nabla_{e_1} e_3 = 0+\phi e_1 = -e_2  \\
  = 0(g(e_1,e_3)e_3 - \eta(e_3)e_1) - 1(g(\phi e_1,e_3)e_3 - \eta(e_3)\phi e_1).\quad \quad \quad \quad
\end{multline}\\
Hence from \eqref{5.1}, \eqref{5.2} and \eqref{5.3} we can see that the manifold $M$ satisfies \eqref{2.5} for $X = e_1, \alpha = 0, \beta = -1,$ and $e_3 = \xi$. Similarly, it can be shown that for $X = e_2$ and $X = e_3$ the manifold also satisfies \eqref{2.5} for $\alpha = 0, \beta = -1,$ and $e_3 = \xi$.\\
Hence the manifold $M$ is a 3-dmensional trans-Sasakian manifold of type (0,-1).
Also, from the definition of the Riemannian curvature tensor $R$ \eqref{1.4}, we get,
$$R(e_1,e_2)e_2 = -e_1, R(e_1,e_3)e_3 = -e_1, R(e_2,e_1)e_1 = -e_2,$$
$$R(e_2,e_3)e_3= -e_2, R(e_3,e_1)e_1 = -e_3, R(e_3,e_2)e_2 = -e_3.$$
Then, the Ricci tensor $S$ is given by,
\begin{equation}\label{5.4}
  S(e_1,e_1) = -2, \quad S(e_2,e_2) = -2, \quad S(e_3,e_3) = -2.
\end{equation}
Then the scalar curvature is,
\begin{equation}\label{5.5}
   r = -6.
\end{equation}
From \eqref{3.5}, we have,
\begin{equation}\label{5.6}
   S(e_1,e_1)=\frac{\lambda+\mu}{2}-(\alpha^2-\beta^2),\smallskip S(e_2,e_2)=\frac{\lambda+\mu}{2}-(\alpha^2-\beta^2),\smallskip S(e_3,e_3)=2(\alpha^2-\beta^2).
\end{equation}
Then from \eqref{5.4} and \eqref{5.6}, we get,
\begin{equation}\label{5.7}
   \frac{\lambda+\mu}{2}-(\alpha^2-\beta^2) = -2,
\end{equation}
and
\begin{equation}\label{5.8}
   \alpha^2-\beta^2 = -1.
\end{equation}
Using \eqref{5.8} in \eqref{5.7}, we obtain,
\begin{equation}\label{5.9}
  \lambda+\mu=-6.
\end{equation}
Then the value of $\lambda+\mu$ in \eqref{5.9} is same as the value of $r$ in \eqref{5.5} and so it satisfies \eqref{3.4}.\\
Hence $g$ defines an $\eta$-Yamabe soliton on a 3-dmensional trans-Sasakian manifold $M$.\\

\end{document}